\documentclass[11pt,reqno,palentino]{amsart}
\usepackage{mypack}
\usepackage{multicol,caption}

               {\endgraf\smallbreak
                \noindent\textbf{Definitie}\sffamily}%
               {\endgraf\smallskip}

  \title{Ropes have an even number of ends}
  \author{Thomas O. Rot}
\begin{document}
\maketitle

\begin{abstract}
  On 13-01-2024 the annual wintersymposium of the Koninlijk Wiskundig Genootschap (KWG) was held in the academiegebouw in Utrecht. The symposium had the theme ``inzichtelijk abstract''. Thomas Rot gave a lecture on his favourite theorem from topology. This article is a written account of this lecture. Audience comprised mostly of high school teachers and that is also the target audience of this article. The slides (in Dutch), which contain more pictures, are available~\cite{T}.
\end{abstract}

\section{My favourite Theorem}
I appreciate simple mathematical ideas which have far-reaching consequences. My favourite theorem in mathematics is such an idea.

\begin{theorem}
Ropes have an even number of ends. 
\label{thm:fav}
\end{theorem}

\begin{figure}
 \centering
 \includegraphics[width=.4\linewidth]{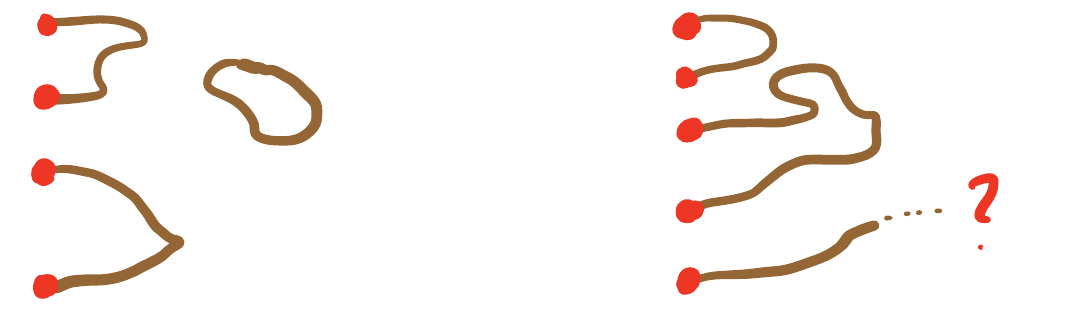}
 \captionof{figure}{An even number of points is the boundary of some ropes, while an odd number of points is not the boundary of some ropes.  }
  \label{fig:even}
\end{figure}

Here and below we assume that ropes have finite length. The proof of this theorem shows that a single rope is either a circle or an interval, which both have an even number of ends. It follows that a finite number of points is the boundary of some ropes if and only if the number of points is even, see Figure~\ref{fig:even}. The formal statement is that the boundary of a compact (finite length) one-dimensional manifold (rope) has an even number of boundary points (ends).

To explain why I love this theorem, I will discuss some of its consequences which at first glance do not seem to have to do anything with ropes. We will see that the theorem can be used to escape a maze, to understand how often subspaces intersect, to explain why the Klein bottle cannot be made to fit into $\mR^3$ without self-intersections and why each closed planar curve contains an inscribed rectangle. I will end with a brief outlook on cobordism theory which is the natural home of Theorem~\ref{thm:fav}.

\section{Escaping a maze}
\begin{figure}
 \centering
 \includegraphics[width=.49\linewidth]{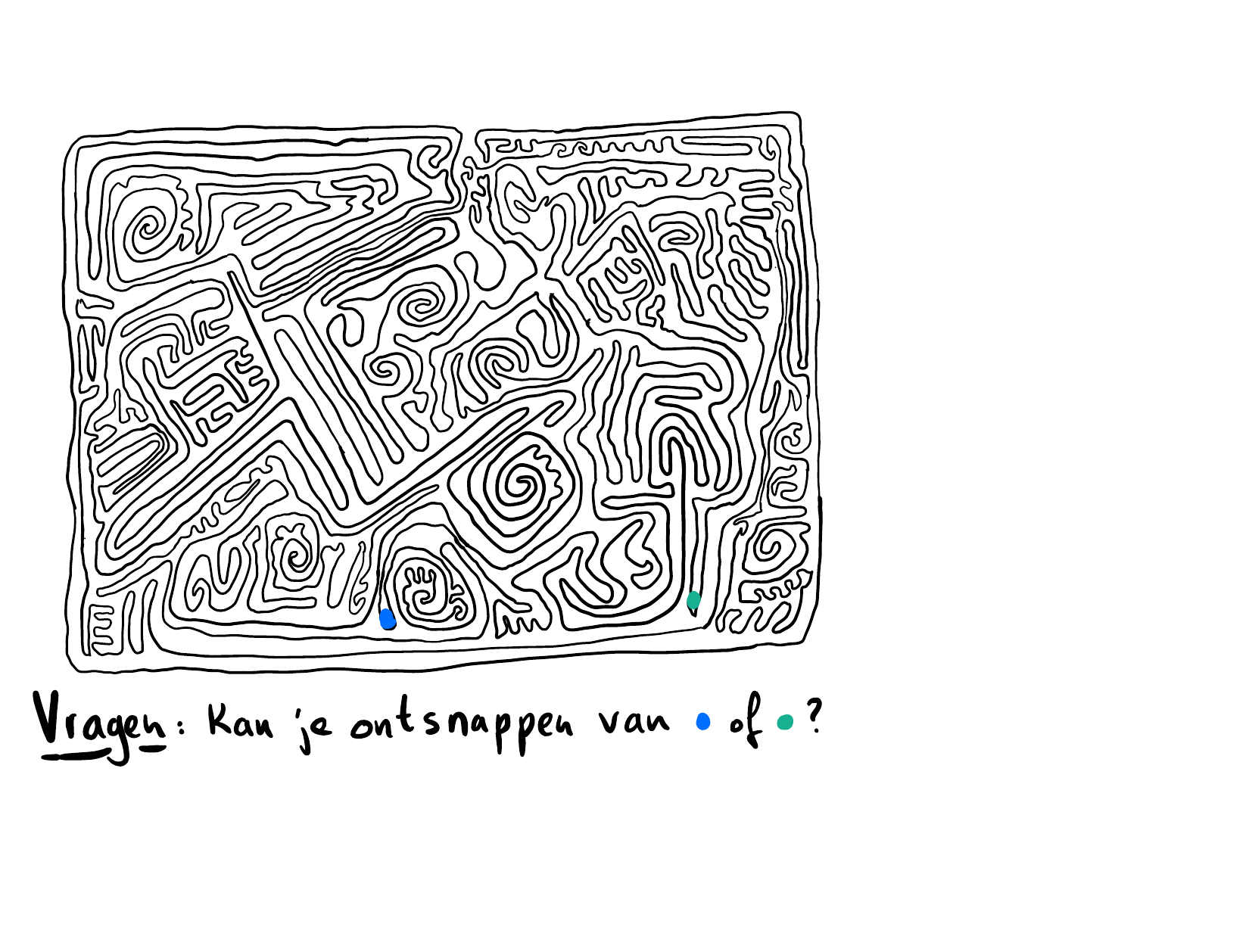}
  \includegraphics[width=.49\linewidth]{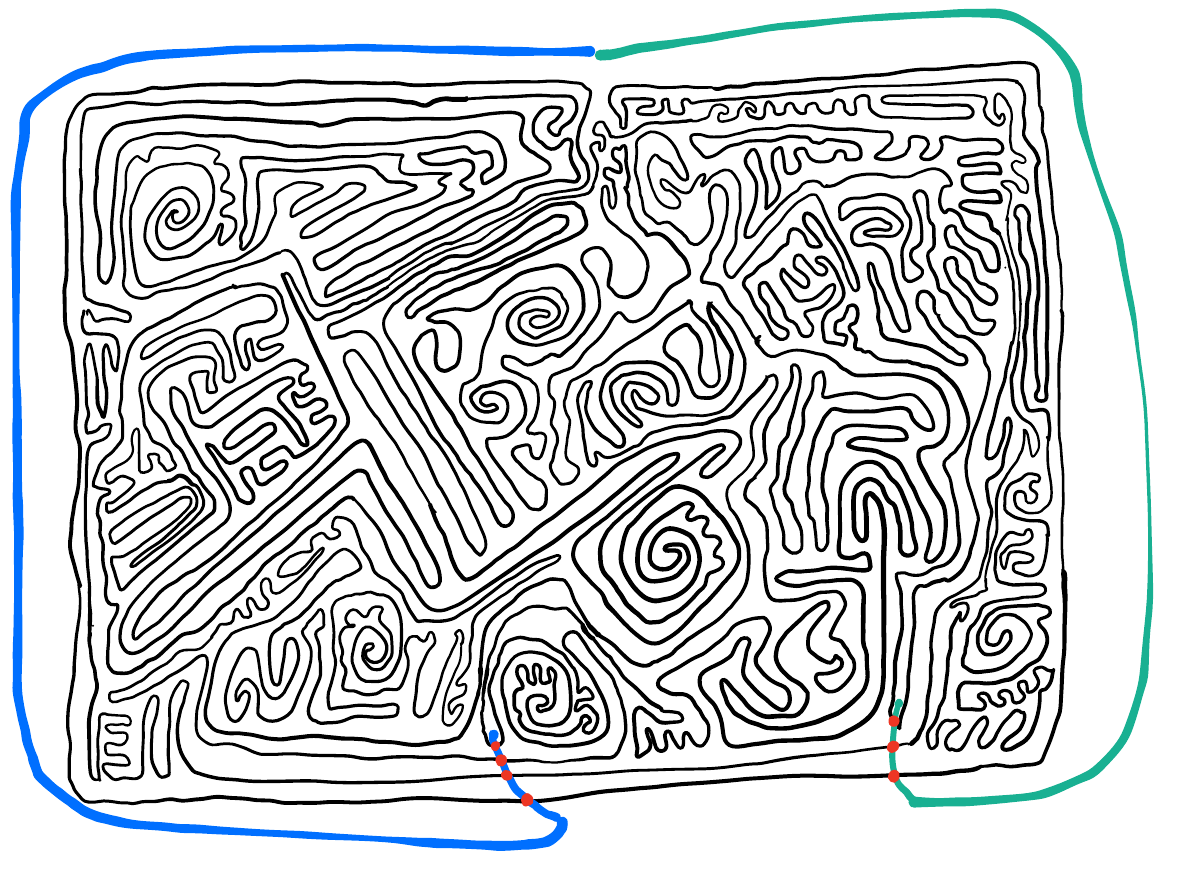}
  \captionof{figure}{How to quickly determine if you can escape a maze as on the left? Draw any curve to the outside and determine if the number of intersections with the maze is even or odd. If it is even, you can escape, if it is odd you can't. In the maze above you can escape from the blue point, as the number of intersections is four (even), but you cannot escape from the turquoise point, as the number of intersections is three (odd). }
  \label{fig:maze}
\end{figure}

\begin{figure}[b]
 \centering
 \includegraphics[width=.3\linewidth]{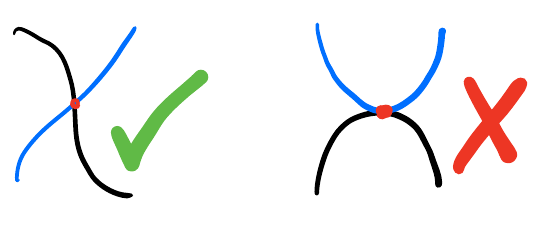}
 \captionof{figure}{Counting intersections should only be done when the curves meet transversely. On the left a transverse intersection is depicted, while on the right a non-transverse intersection is depicted: the two curves are tangent to each other. A somewhat deep fact in differential topology is that generically curves only have transverse intersections. By an arbitrary small perturbation of the curve all intersections will be transverse, and transversality is preserved under small enough perturbations.}
  \label{fig:transverse}
\end{figure}

\begin{figure}
 \centering
 \includegraphics[width=.7\linewidth]{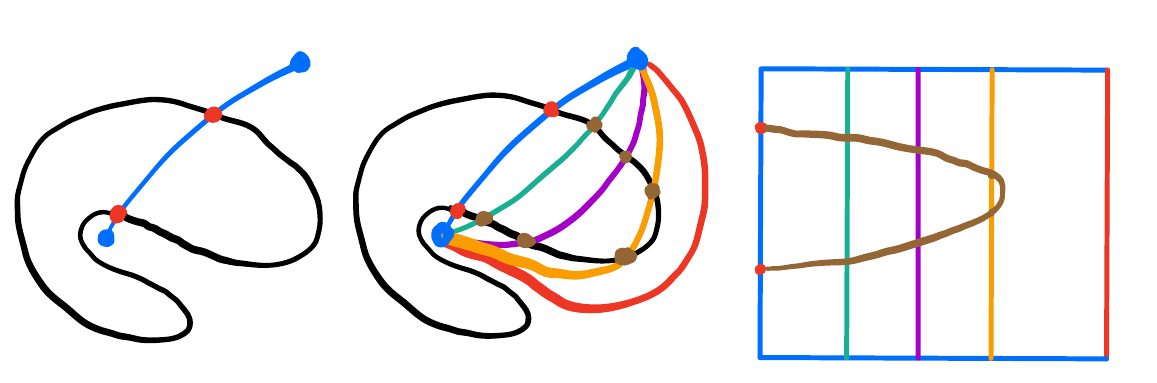}
 \captionof{figure}{Can we escape this very simple maze on the left? Yes: the curve in blue intersects the maze twice, which is even. Why is this number even? Imagine deforming the blue curve to the red curve as in the middle picture. The red curve is an actual escape route. We can keep track of the intersections with the maze during this deformation. This is depicted on the right: In brown the intersections of all the curves with the maze are drawn. These are ropes, and the ends of these ropes are the intersections of the original curve with the maze. This shows that the intersection is even if we can escape. A slightly more complicated argument shows that this condition is also sufficient. To avoid a flood of complaints I need to be precise: we assume that the maze is a connected closed simple curve.}
  \label{fig:homotopy1}
\end{figure}

How can you determine quickly if you can escape a maze such as in Figure \ref{fig:maze}? There is a trick with which you can answer this question very quickly: Draw any\footnote{Any \emph{generic} curve to be precise, see Figure~\ref{fig:transverse}.} curve from the starting point somewhere in the maze to the outside. Count the number of intersections of this curve with the maze. If the number of intersections is even, you can escape, but if the number of intersections is odd, you cannot escape. Why this works is explained in Figure~\ref{fig:homotopy1}.

What we have shown is that the parity the number of intersections between a curve with a closed curve is constant under deformations which keep the endpoint of the curve fixed. A similar argument also shows that any two closed curves in the plane intersect each other in an even number of points. In turn this statement can be generalized to higher dimensions. Here is the statement in three dimensions.

\begin{theorem}
  \label{thm:three}
A closed curve and a closed surface intersect in an even number of points in three dimensions. 
\end{theorem}

\begin{figure}
 \centering
 \includegraphics[width=\linewidth]{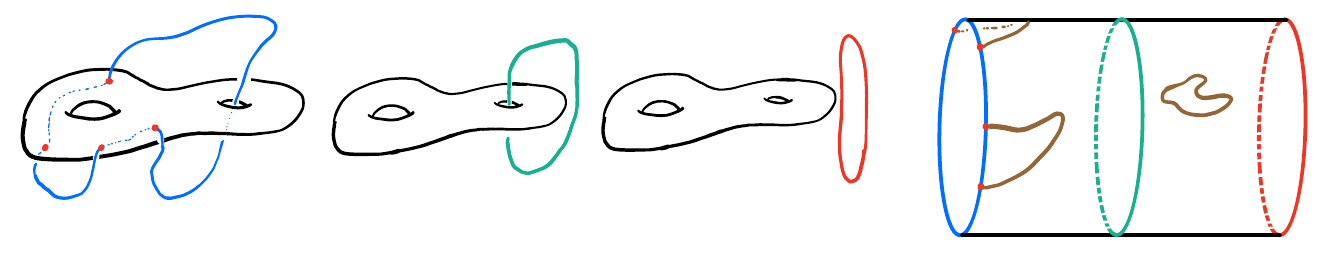}
 \captionof{figure}{The number of intersections of a closed surface and a closed curve in three dimensions is even. The proof is analogous of the proof described in Figure~\ref{fig:homotopy1}: deform the curve to lie outside the surface and keep track where the family of curves intersect the surface. This traces out ropes, whose ends are the original intersections, which must therefore be even. Note that one of the ropes is a closed curve itself, which did not occur in Figure~\ref{fig:homotopy1}. }
  \label{fig:threedim}
\end{figure}

A surface is called closed if it does not have a boundary and if it is compact. Compactness intuitively means that the surface does not extend out to infinity. In Figure~\ref{fig:threedim} a proof of Theorem~\ref{thm:three} is sketched.

\section{Non-orientable surfaces}

\begin{figure}[b]
 \centering
 \includegraphics[width=\linewidth]{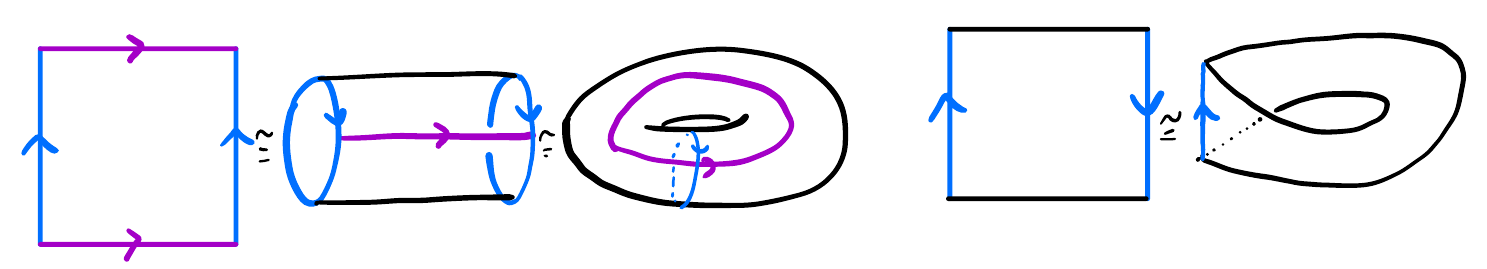}
 \captionof{figure}{Gluing edges of a polygon pairwise is a way to construct surfaces. On the left the torus is constructed, which is a closed surface (without boundary). The M\"obius strip, does have a boundary. It is also not orientable, cf.~\ref{fig:escher}.}
  \label{fig:mobius}
\end{figure}

\begin{figure}[t]
 \centering
 \includegraphics[width=.3\linewidth]{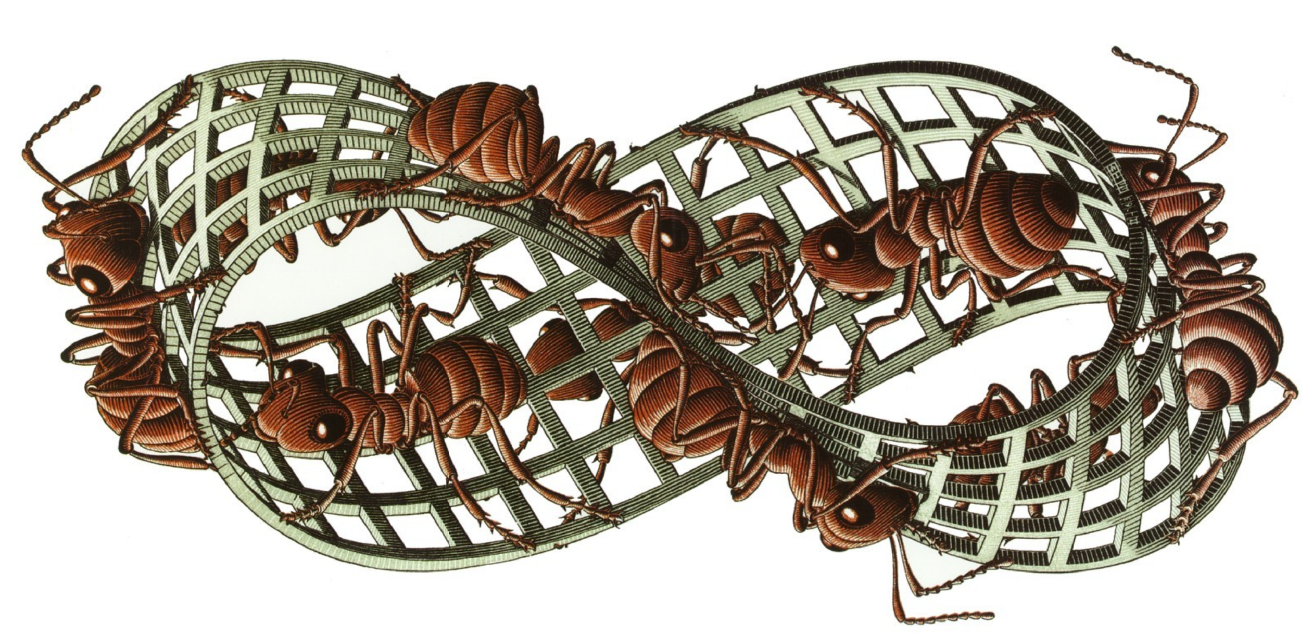}
 \captionof{figure}{The M\"obius strip is non-orientable. Taking a roundtrip around the M\"obius strip makes the ants get to the other side of the strip. This is a famous work of Escher~\cite{Escher}.}
  \label{fig:escher}
\end{figure}

Compact surfaces can be constructed by a gluing process: Take a finite number of polygons, and glue the edges pairwise. If each edge is glued to another one a closed surface is constructed, if not all edges are glued a compact surface with boundary is made. An interesting example of a compact surface with boundary is the M\"obius strip, see Figures~\ref{fig:mobius} and~\ref{fig:escher}.

\begin{figure}[b]
 \centering
 \includegraphics[width=\linewidth]{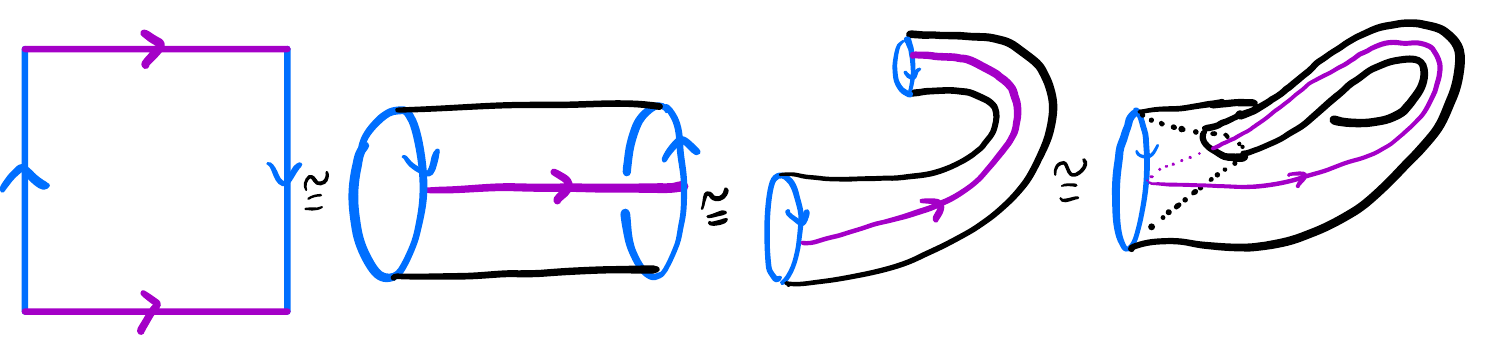}
 \captionof{figure}{Identifying the square as on the left produces the Klein Bottle. The Klein bottle is a closed non-orientable surface and cannot be put inside three dimensions. The Klein bottle contains a M\"obius strip: The complement of a small neighborhood of the purple curve is a M\"obius strip}
  \label{fig:klein}
\end{figure}

The M\"obius strip is not orientable and has a boundary. The Klein bottle is an example of a \emph{closed} non-orientable surface, see Figure~\ref{fig:klein}. Depicting the Klein bottle is harder than depicting the M\"obius strip. It is not possible to fit, or \emph{embed}, the Klein bottle in three dimensions without self intersections. Let me state this as a theorem. 

\begin{theorem}
The Klein bottle cannot be embedded in $\mR^3$. 
\end{theorem}

\begin{figure}
 \centering
 \includegraphics[width=.8\linewidth]{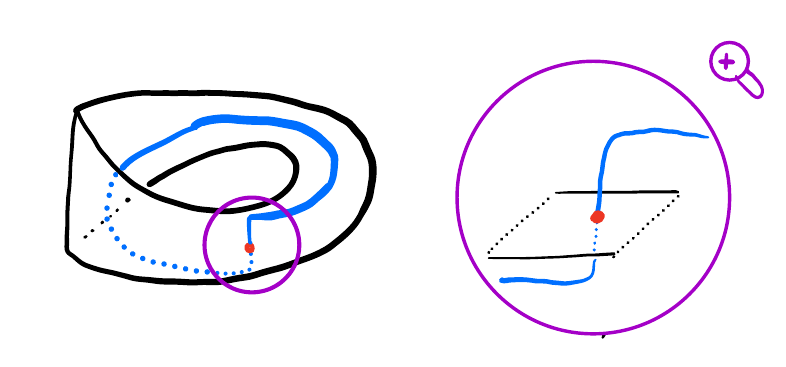}
 \captionof{figure}{In a neighborhood of a M\"obius strip in three dimensions we can find a loop that intersects the M\"obius strip once. This does not contradict Theorem~\ref{thm:three} as the M\"obius strip has a boundary and is therefore not closed. The Klein bottle contains a M\"obius strip. If the Klein bottle would embed in three dimensions, we can construct a blue curve following the M\"obius strip contained in the embedded Klein bottle. The blue curve can be made not to intersect the Klein bottle in any other point, as it can be made arbitrary close to the M\"obius strip. Thus the number of intersections of the Klein bottle and the curve is odd, which contradicts Theorem~\ref{thm:three}. The Klein bottle therefore cannot be embedded in three dimensions. }
  \label{fig:samelson}
\end{figure}

We will follow the wonderful argument of Samelson~\cite{S} to prove this. Imagine, with a contradiction in mind, that the Klein Bottle \emph{does} embed in three dimensions. Then in  Figure~\ref{fig:samelson} it is shown that if the Klein bottle would embed, that there would also have to exist a closed curve which intersects the Klein bottle \emph{once}. But Theorem~\ref{thm:three} tells us that this is impossible. Hence the Klein bottle does not embed in $\mR^3$. Note that the existence of the curve that intersects the M\"obius strip once does not contradict Theorem~\ref{thm:three} as the M\"obius strip has a boundary and is therefore not closed.

Any surface that is non-orientable contains a M\"obius strip. The argument just given proves the following Theorem.

\begin{theorem}
  \label{thm:nonorientable}
A closed non-orientable surface does not embed in $\mR^3$. 
\end{theorem}

Even though it is possible to embed the M\"obius strip in three dimensions, we cannot prescribe the behavior of the boundary of the M\"obius strip completely at will. For example we have the following corollary, which we will need later.

\begin{corollary}
  \label{cor:mobius}
The M\"obius strip cannot be embedded in three dimensions in such a way that the boundary of the M\"obius strip is contained in the plane $\{(x,y,z)\,|\,z=0\}$ and the interior of the M\"obius strip is contained in the upper half space $\{(x,y,z)\in \mR^3\,|\,z>0\}$.
\end{corollary}

To see that this corollary follows from Theorem~\ref{thm:nonorientable} one imagines capping off the curve with a disc in the lower half-space $\{(x,y,z)\in \mR^3\,|\,z<0\}$ in a smooth manner. This produces an embedded non-orientable closed surface in three dimensions which is not possible by Theorem~\ref{thm:nonorientable}. Thus the M\"obius strip cannot be embedded in this way. 

\section{A surprising M\"obius strip}

\label{sec:surprise}
\begin{figure}
  \begin{minipage}{.55\textwidth}
 \includegraphics[width=.8\linewidth]{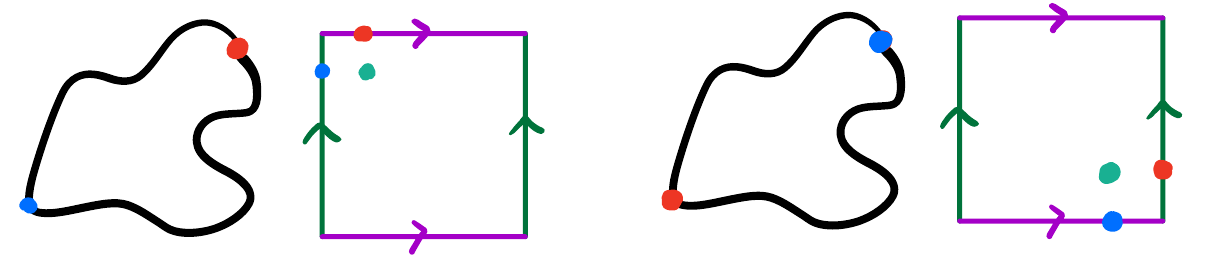}
 \captionof{figure}{The space of \emph{ordered} points on the circle is a torus. The $x$ and $y$ coordinates of a point on the square in the picture each represent a point on the curve. The torus is obtained by gluing the square as prescribed. In the picture the configuration of the points on the left is different from the configuration of the points on the right as we can distinguish which point is red, and which point is blue, see also Figure~\ref{fig:unordered1}.}
  \label{fig:ordered}
\end{minipage}%
\begin{minipage}{.44\textwidth}
 \includegraphics[width=\linewidth]{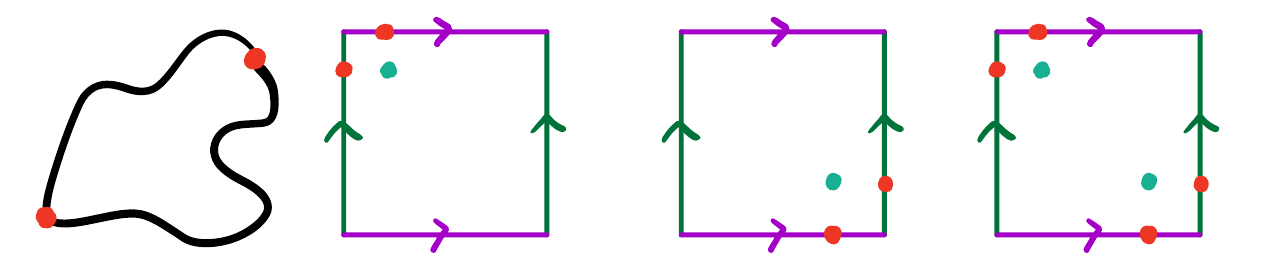}
 \captionof{figure}{The space of \emph{unordered} pairs of points on the circle can be visualized as the space of pairs of points on the torus, where the points are identified if they are mirrored along the diagonal. }
  \label{fig:unordered1}
\end{minipage}
\end{figure}

A critique topologists often have to answer to is that they play with toys that they invent themselves. \emph{It is nice that non-orientable surfaces exist, but as they cannot be found in three dimensions, what is the point of studying them?} Here I would like to argue that the M\"obius strip is not only a mathematical curiosity for the sake of it, but naturally occurs in the world around us. Imagine a particle constrained on a closed curve in the plane. Mathematically we can parametrize this particle by a circle, and a particle is then just a point on the circle. What about two particles? If we can distinguish the particles from each other the particles can be parametrized by the product of two circles, which is a torus, see Figure~\ref{fig:ordered}. But if the particles are indistinguishable this is not correct, see Figure~\ref{fig:unordered1}. We need to make identifications on the torus, as multiple points on the torus correspond to the same configuration if we cannot distinguish the particles. If we identify these points we obtain a M\"obius strip! In Figure~\ref{fig:unordered2} a graphical proof of this fact is shown. The M\"obius strip occurs naturally as the configuration space of pairs of unordered points on a curve! We will now use this fact in an unexpected way. 



\section{Pegs in curves}

We now turn to discuss the proof, due to Herbert Vaughan, of a curious fact in plane geometry. The Nieuw Archief has already paid attention to the history of this result: Robbert Fokkink gave a beautiful historical account~\cite{Fok}. I do not have anything to add to the historical account, but I do want to give a few more details on the proof as it fits nicely with the mathematics we have discussed so far.

A priori this theorem does not have anything to do with ropes and non-orientable surfaces. The statement that we will prove is the following:

\begin{theorem}
  \label{thm:inscribed}
Every smooth simple closed curve in the plane contains four distinct vertices that form a non-degenerate rectangle. 
\end{theorem}

\begin{figure}[b]
 \centering
 \includegraphics[width=.3\linewidth]{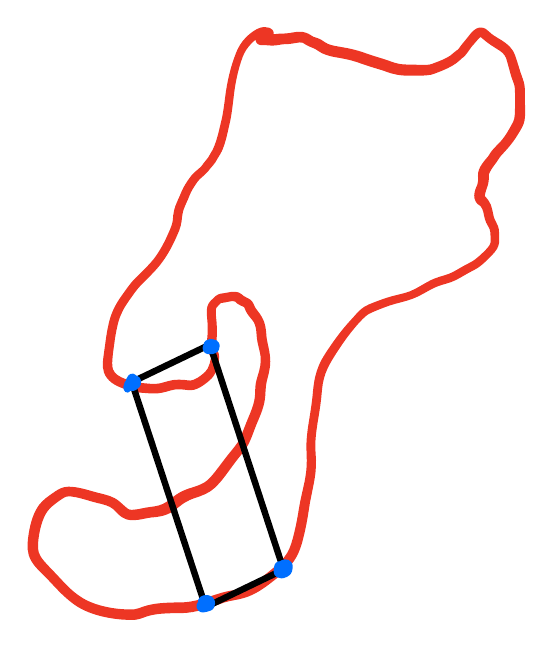}
 \captionof{figure}{Greene and Lobb proved that every smooth simple closed curve in the plane contains uncountably many rectangles whose vertices lie on the given curve. We prove the existence of only one.}
  \label{fig:rectangle}
\end{figure}

See Figure~\ref{fig:rectangle} for an example of such a rectangle. We call these rectangles \emph{inscribed} rectangles. Recently there is renewed interest in this problem as Greene and Lobb \cite{GL} proved a much stronger result:
\begin{figure}
 \centering
 \includegraphics[width=\linewidth]{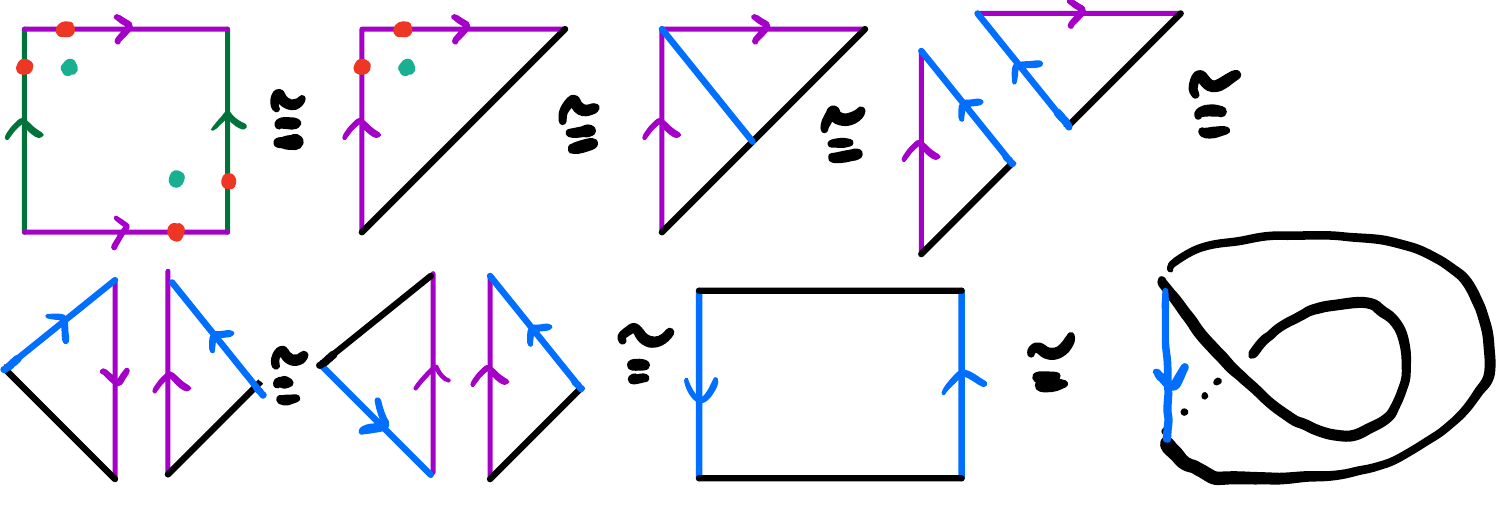}
 \captionof{figure}{The space of \emph{unordered} pairs of points on the circle is the M\"obius strip. By Figure~\ref{fig:unordered1} we can see this space as the space of two points on the torus ``mirrored'' on the diagonal. To parametrize these points we only need to keep track of the point in the upper triangle. We do need to keep in mind that the green line now will be identified with the purple line. Cutting the obtained surface along the blue line, and gluing it back reveals that this space was the M\"obius strip in disguise.}
  \label{fig:unordered2}
\end{figure}

\begin{theorem}
Every smooth simple closed curve contains uncountably many inscribed rectangles, at least one for each aspect ratio of the long and short side of the rectangle. 
\end{theorem}

The proof of this theorem is outside the scope of this article, but see~\cite{Fok} for a sketch of the proof\footnote{I should remark about a small typo in the end of the article~\cite{Fok}. For the contradiction Greene and Lobb construct a \emph{Lagrangian} Klein bottle, not a symplectic one.}.

\begin{figure}[b]
 \centering
 \includegraphics[width=.3\linewidth]{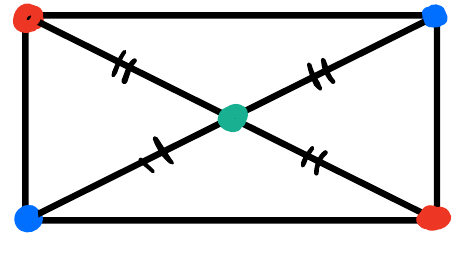}
 \captionof{figure}{Two pairs of points (the red and blue points) form a rectangle if and only if the midpoint of the red points equals the midpoint of the blue points and the distance to the common midpoint (turquoise) is the same. This is an exercise in Euclidean geometry. }
  \label{fig:euclidean}
\end{figure}

We will prove Theorem~\ref{thm:inscribed} now. The only ingredient missing is a fact from Euclidean Geometry. To determine if four points form a rectangle, it suffices to show that pairs of these points have a common midpoint and the distance from all points to the midpoint is the same, see Figure~\ref{fig:euclidean}. Let $O$ be the space of pairs of points on the closed curve. Thus $(x,y)\in O$ means that $x\in \mR^2$, $y\in\mR^2$ and $x,y$ both lie on the curve. This space $O$ is topologically the same as a torus by the discussion in Section~\ref{sec:surprise}.

\begin{figure}
 \centering
 \includegraphics[width=\linewidth]{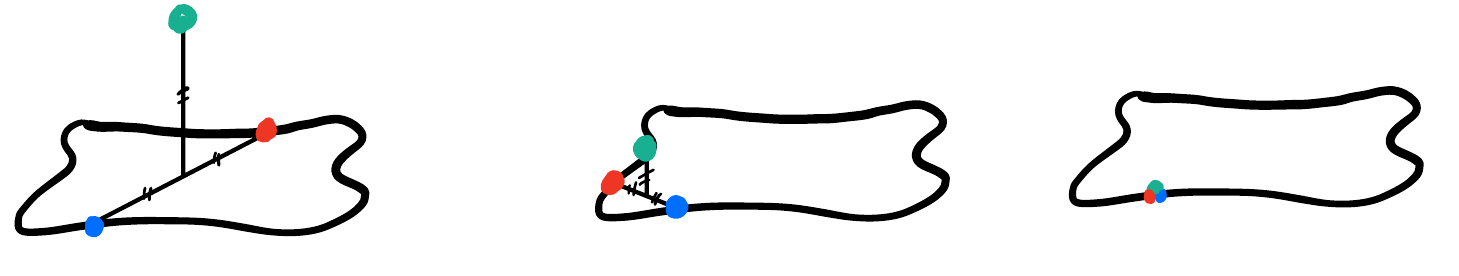}
 \captionof{figure}{Some values of the function $f:O\rightarrow \mR^3$. The red and blue point are mapped to the point above the midpoint, with height the distance to the midpoint. Note that this function is invariant under swapping the red and blue point, and if the red and blue point get closer to each other, the height goes to zero. }
  \label{fig:function}
\end{figure}

Define the function $f:O\rightarrow \mR^3$ by
\begin{equation}
  f(x,y)=\left(\frac12(x+y),\frac{1}{2}\norm{x-y}\right).
\end{equation}
The expression $\frac12(x+y)\in \mR^2$ is the midpoint between the two points $x$ and $y$ on the curve, and $\frac{1}{2}\norm{x-y}$ is the distance from the points to the midpoint. In Figure~\ref{fig:function} some values of the function are depicted. A nice movie about this function can be found on the YouTube channel of Grant Sanderson~\cite{3blue1brown}. The function has three important properties: i) It satisfies $f(x,y)=f(y,x)$, so we can view $f:U\rightarrow \mR$ as a function on the space $U$ of \emph{unordered} points on the curve. Recall that this space is topologically a M\"obius strip! ii) The height of the point $f(x,y)$ is always $\geq 0$ with equality if and only if $x=y$. iii) If $f(x,y)=f(x',y')$ then the points $x,x',y,y'$ lie on a rectangle, see Figure~\ref{fig:euclidean}. \\

Summarizing: The map $f$ maps the M\"obius strip $U$ to $\mR^3$, where the boundary of the strip is mapped to the curve in the $x-y$ plane and the interior of the M\"obius strip has postive $z$ coordinate. If there is \emph{no} inscribed rectangle on the curve, the map $f$ is injective and gives an embedding of the M\"obius strip in a way that contradicts Corollary~\ref{cor:mobius}. As this contradicts the Corollary, the hypothesis that $f$ is injective is false. Thus an inscribed rectangle must exist!


\section{The cobordism ring and a theorem of Thom}

I would like to end by putting my favourite Theorem~\ref{thm:fav} in a broader mathematical context. Curves and surfaces are examples of smooth manifolds: spaces that locally resemble Euclidean space $\mR^n$, and which have a notion of differentiability. We also have manifolds with boundary, which locally resemble $\mR^n\times[0,\infty)$. Closed manifolds, manifolds which are compact and do not have a boundary, are of particular interest to topologists.

Classifying closed manifolds up to diffeomorphism, the natural notion of sameness of manifolds, is very hard. In low dimensions we can make progress: A closed zero-dimensional manifold is a finite number of points, and two zero-dimensional manifolds are diffeomorphic if and only if they have the same number of points. A connected closed one-dimensional manifold is diffeomorphic to a circle, and one-dimensional closed manifolds are classified by the number of connected components.

\begin{figure}
 \centering
 \includegraphics[width=.8\linewidth]{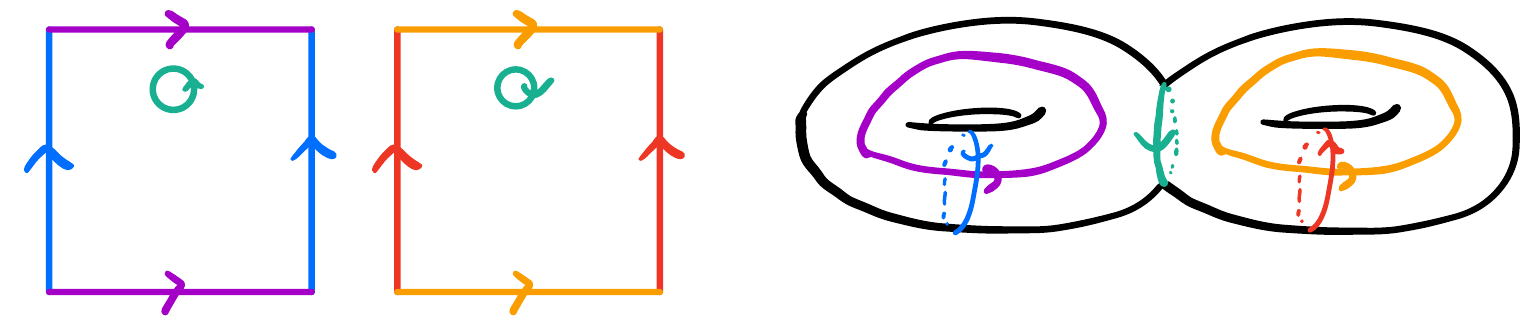}
 \captionof{figure}{The connected sum of two tori is a genus-two surface. All orientable surfaces are either the two-sphere, the torus or connected sums of tori. All non-orientable surfaces are either the real projective plane $\mR\mP^2$, or connected sums of those.}
  \label{fig:connectedsum}
\end{figure}

Connected two-dimensional surfaces come in two families. The first family are the orientable surfaces. These are the two dimensional sphere, the torus, and the genus-$g$ surfaces.  The connected sum of two surfaces is obtained by removing two small discs and gluing the resulting surfaces along the new boundaries together. The genus-$g$ surfaces can be obtained by repeatedly taking a connected sum with a torus, see Figure~\ref{fig:connectedsum}. There is also a family of non-orientable surfaces. These are constructed by connected sums of the real projective plane $\mR\mP^2$. The real projective plane is obtained by gluing a square as in Figure~\ref{fig:klein} but flipping the orientation of one of the purple edges. The Klein bottle is the connected sum of $\mR\mP^2$ with itself. All non-orientable surfaces can be constructed by taking further connected sums with $\mR\mP^2$. 

In three dimensions one also can make progress in the classification, but this is much harder. After three, the music stops. The diversity of four-dimensional manifolds is staggeringly vast. It is so big that one can show that in a precise sense \emph{no} reasonable classification of closed four-dimensional manifolds is possible.

One way to proceed with a classification of manifolds is to take a coarser notion of sameness. If more things are declared the same, fewer distinct classes remain, and it might be possible to classify them. For example we can classify manifolds by their dimension, but this relation is too coarse to be really useful in the study of manifolds. Thom~\cite{Thom} surprised the mathematical community in the 50-ies by classifying closed manifolds in every dimension up to \emph{cobordism}. Two closed manifolds $M,N$ are cobordant if there exists a compact manifold with boundary $W$ such that the boundary $\partial W$ is the disjoint union of $M$ and $N$, see Figure~\ref{fig:cobordism}. Diffeomorphic manifolds are cobordant, so this is indeed a coarser notion of sameness compared to diffeomorphism. Using our classification of low-dimensional manifolds, we also get a classification of manifolds up to cobordism in low dimensions.

The following is a reformulation of Theorem~\ref{thm:fav} from the beginning of this article.
\begin{theorem}
 Every closed zero-dimensional manifold is either cobordant with the empty manifold, or with the manifold with one point. There are two cobordism classes of zero-dimensional manifolds: the class of ``even number of points'' and the class of ``odd number of points''.
\end{theorem}

Dimension one is boring. Every circle is the boundary of a disc. Thus the circle is cobordant with the empty manifold. We say that the circle is nullbordant. As any closed one dimensional manifold is a finite number of circles we get the following Theorem. 

\begin{theorem}
Every one-dimensional manifold is nullbordant. There is only one cobordism class of one-dimensional manifolds.
\end{theorem}

Dimension two is interesting. All orientable surfaces are nullbordant. To see this one can embed the surface in $\mR^3$ and take the inside as a bounding three dimensional manifold.

Half of the non-orientable surfaces are nullbordant. The Klein Bottle for example \emph{is} the boundary of a three-dimensional manifold. Imagine the Klein bottle as a family of circles parametrized by a circle as in the middle picture of Figure~\ref{fig:klein}. Filling in these circles with discs produces a three-dimensional compact manifold whose boundary is the Klein bottle.\\

The surface $\mR\mP^2$ is not nullbordant, which is harder to see. It is possible to use my favourite Theorem~\ref{thm:fav} to prove this (for the initiated the Euler characteristic modulo $2$ is a cobordism invariant), but it would take us too far astray to discuss this further. All surfaces are either cobordant with the Klein bottle, or with $\mR\mP^2$, which gives:

\begin{figure}
 \centering
 \includegraphics[width=.4\linewidth]{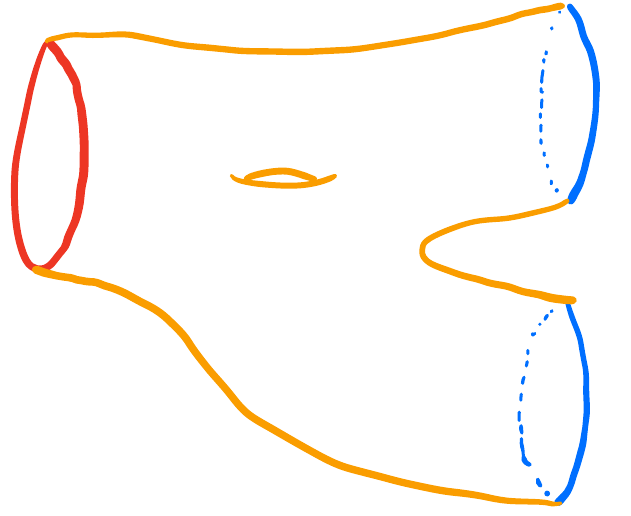}
 \captionof{figure}{The red manifold and the blue manifold (which has two components) are cobordant. There exists a manifold, of one dimension higher, whose boundary is the union of the red and blue manifold. Another example seen on the left of Figure~\ref{fig:even}. The zero-dimensional manifold with an even number of points is nullbordant (cobordant with the empty set). On the right it is shown that an odd number of points is \emph{not} nullbordant.}
  \label{fig:cobordism}
\end{figure}

\begin{theorem}
  There are two cobordism classes of two-dimensional manifolds. Every closed two-dimensional manifold is either nullbordant, or cobordant with $\mR\mP^2$. 
\end{theorem}

In these three theorems on the cobordism classes we heavily used the classification of the manifolds \emph{up to diffeomorphism}. It was Thom's great insight~\cite{Thom} that manifolds of \emph{all} dimensions can be classified up to cobordism, \emph{even though we cannot classify the manifolds themselves}!

The set of cobordism classes can be equipped with a rich algebraic structure. We can define an addition operation by the disjoint union of manifolds, and a multiplicative operation by the cartesian product of manifolds. These operations satisfy properties like commutativity, associativity and distributivity, just like ordinary addition and multiplication. This turns the set $\mathcal{N}_*$ of cobordism classes of closed manifolds into a (graded) ring\footnote{We listen to the advice of Beyonc\'e~\cite{be}. We like it, so we put a ring on it.}. The unit $0$ for addition is represented by the empty manifold, and the unit $1$ for multiplication is represented by a point. An interesting feature of the cobordism ring is that  $[M]+[M]=0$ for all cobordism classes $[M]$, as two disjoint copies of a closed manifold $M$ is the boundary of the ``cylinder'' $M\times[0,1]$. This is similar to calculating modulo $2$, i.e.~in $\mathbb{Z}/2\mathbb{Z}$. Thom's theorem is now a computation of the ring of cobordism classes of manifolds. 

\begin{theorem}
 The cobordism ring $\mathcal{N}_*$ is isomorphic to the polynomial algebra $
\mathbb{Z}/2\mathbb{Z}[x_2,x_4,x_5,x_6,x_8,\ldots]
$
with one generator $x_k$ in dimension $k$ for each $k\not=2^i-1.$
\end{theorem}

Much more is known about the cobordism ring and variants of it, but many open questions remain. Cobordism theory remains central in algebraic and differential topology. 

\section{Acknowledgements}

I would like to thank Renee Hoekzema, Michael Jung and Mark Timmer for comments on a draft of this article. I would also like to thank Carmen Oliver Huidobro for discussions during the writing of her bachelor thesis on Greene and Lobb's theorem.


\end{document}